\documentclass{ifacconf}

\usepackage{graphicx}      
\usepackage{natbib}        

\usepackage{amsmath}
\usepackage{amsfonts}
\usepackage{amssymb}
\usepackage{color}
\usepackage{dsfont}
\usepackage{mathrsfs}
\usepackage[usenames,dvipsnames]{pstricks}
\usepackage{enumerate} 
\usepackage{easyReview} 

\allowdisplaybreaks
\usepackage{flushend}

\usepackage{comment}
\usepackage{color} 
\usepackage{subcaption}

\begin{document}
\begin{frontmatter}

\title{Neural Operators for Adaptive Control of Traffic Flow Models\thanksref{footnoteinfo}} 

\thanks[footnoteinfo]{This work was supported by the National Natural Science Foundation of China under Grants 92471108, 12131008.}

\author[First]{Kaijing Lyu} 
\author[First]{Junmin Wang} 
\author[Fourth]{Yihuai Zhang} 
\author[Fourth]{Huan Yu} 
\address[First]{School of Mathematics and Statistics,  Beijing Institute of Technology, Beijing,China, 100081(e-mail: kjlv@bit.edu.cn)}
\address[Fourth]{Hong Kong University of Science and Technology (Guangzhou), Nansha, Guangzhou, Guangdong, China, 511400(e-mail: huanyu@ust.hk)}

\begin{abstract}                
The uncertainty in human driving behaviors leads to stop-and-go instabilities in freeway traffic. The traffic dynamics are typically modeled by the Aw-Rascle-Zhang (ARZ) Partial Differential Equation (PDE) models, in which the relaxation time parameter is usually unknown or hard to calibrate. This paper proposes an adaptive boundary control design based on neural operators (NO) for the ARZ PDE systems.  
In adaptive control, solving the backstepping kernel PDEs online requires significant computational resources at each timestep to update estimates of the unknown system parameters.
To address this, we employ DeepONet to efficiently map model parameters to kernel functions. Simulations show that DeepONet generates kernel solutions nearly two orders of magnitude faster than traditional solvers while maintaining a loss on the order of \(10^{-2}\). Lyapunov analysis further validates the stability of the system when using DeepONet-approximated kernels in the adaptive controller. This result suggests that neural operators can significantly accelerate the acquisition of adaptive controllers for traffic control.
\end{abstract}


\end{frontmatter}

\section{Introduction}
Stop-and-go traffic congestion causes unsafe driving, pollution, and delays \citep{belletti2015prediction}, driven by the propagation of shock waves from delayed driver responses. The ARZ model effectively describes this dynamic \citep{aw2000resurrection,zhang2002non}. 
PDE backstepping has been applied to traffic control using state and output feedback \citep{yu2019traffic}. Recent studies extend it to multilane and mixed autonomy traffic \citep{burkhardt2021stop,yu2018adaptive,yu2022traffic}. 
In traffic flow modeling, relaxation time represents drivers' reaction delays, but its heterogeneity and unpredictability make it practically unmeasurable. The adaptive control method has been well studied to stabilize PDEs with unknown parameters\citep{krstic2008adaptive,anfinsen2019adaptive,yu2018adaptive,yu2018varying}. However, implementing adaptive control for the ARZ model remains challenging due to high computational demands for real-time parameter and state estimation. This is because the adaptive control process simultaneously requires the estimation of unknown system parameters and PDE states. After each time step, it is necessary to recalculate the solution of the gain kernel functions, which places extremely high demands on real-time computation.To address this, we use neural operators to accelerate the computation of adaptive controller of ARZ PDE system.

\par Recent studies have successfully applied DeepONet to one-dimensional transport PDEs \citep{Bhan23}, reaction-diffusion equations, and observer designs \citep{krstic2024neural}, demonstrating the stability of the system with approximate kernels. Further developments focus on hyperbolic PDEs \citep{qi2024neural}, parabolic PDEs with delays \citep{WangSSHeat23}, and the ARZ PDE system for traffic control \citep{zhang2024neural}. DeepONet is used for gain scheduling in real-time control of nonlinear PDEs \citep{lamarque2024gain}, integrating PDE backstepping for offline learning and improving closed-loop stability in complex PDEs.  
The use of NO-approximated gain kernels for adaptive control was first explored for first-order hyperbolic PDEs \citep{lamarque2024adaptive} and extended to reaction-diffusion equations \citep{bhan2024adaptive2}. 
\par The main contributions are summarized as follows: We propose an adaptive control method based on NO to stabilize the ARZ traffic PDE model with unknown relaxation time.
To address computational challenges in solving gain kernel equations, we integrate DeepONet into the adaptive control framework.  
DeepONet is shown to be nearly two orders of magnitude faster than traditional PDE solvers for generating kernel functions, with a loss of around \(10^{-2}\).  Additionally, we prove system stability through Lyapunov analysis when using DeepONet-approximated kernels in the adaptive controller. This is the first study to combine DeepONet with adaptive control in traffic flow systems, demonstrating its potential to improve computational efficiency in traffic control.

\section{{Nominal Adaptive Control Design for ARZ Traffic }}\label{sec_Nominal}
The ARZ PDE model is used to describe the formation and dynamics of the traffic oscillations which refer to variations of traffic density and speed around equlibrium values. The ARZ model is given by:
\begin{align} 
	\partial_{t} \rho_1+\partial_{x}(\rho_1 v_1) & =0, \label{eq:ARZ_1} \\
	\partial_{t} (v_1-V(\rho_1))+v\partial_{x}\left(v_1-V(\rho_1)\right)  & =\frac{V(\rho_1)-v_1}{\tau} ,\\
	\rho_1(0, t) & =\frac{q_1^{\star}}{v_1(0, t)}, \\
	v_1(L, t) & =U(t)+v_1^{\star},  \label{eq:ARZ_4} 
\end{align}
where $(x,t)\in [0,L]\times \mathbb{R}_+$, 
$\rho_1(x, t)$ represents the traffic density, $v_1(x, t)$ represents the traffic speed. The relaxation time $\tau$ describes how fast drivers adapt their speed to equilibrium speed-density relations. Its value is usually difficult to measure in practice and is easily affected by various external factors. The variable $p_1(\rho_1)$, defined as the traffic system pressure, is related to the density by the equation $p_1(\rho_1)=c_0(\rho_1)^{\gamma_0},$
and $c_0,\gamma_0 \in \mathbb{R}_+$. The equilibrium velocity-density relationship $V (\rho_1)$ is given in Greenshield model:
\begin{equation}
	V(\rho_1)=v_{f}\left(1-\left(\frac{\rho_1}{\rho_{m}}\right)^\gamma_0\right),
\end{equation}
where $v_f$ is the free flow velocity, $\rho_m$ is the maximum density, and $(\rho_1^*,v_1^*)$ are the equilibrium points of the system with $v_1^*=V(\rho_1^*)$. We consider a constant traffic flux $q_1^*=\rho_1^* V(\rho_1^*)$ entering the domain from $x = 0$. We can apply the change of coordinates to rewrite it in the Riemann coordinates and then map it to a decoupled first-order 2$\times$2 hyperbolic system 
\begin{align}\label{eq:uv_def}
	\partial_{t} u(x, t)  & =-\lambda \partial_{x} u(x, t), \\
	\partial_{t} v(x, t)& =\mu \partial_{x} v(x, t) +c(x) u(x, t), \label{eq:uv_def_2}\\
	u(0, t) & =r v(0, t), \\
	v(1, t) & =U(t), \label{eq:uv_def_v0}
\end{align}
where $\lambda=v_1^*$, $\mu=\gamma_0 p_1^*-v_1^*$, ${c}(x)=-\dfrac{1}{\tau}\exp(-\dfrac{x}{\tau v_1^*})$, $r=\frac{\rho_1^*V(\rho_1^*)+v_1^*}{v_1^*}$. The relaxation time $\tau$ is unknown, so the spatially variable coefficient $c(x)$ is also unknown. We assume $c(x)\in C^1([0,1])$ is bounded with $c(x)\leq \bar{c}$. Then, we propose the adaptive control law using passive identifier method. 

We consider the identifier
\begin{align}\label{eq:uv_identifier_u}
	\partial_{t} \hat{u}(x, t) =&-\lambda \partial_{x}\hat{u}(x, t)+\rho e(x, t)\|\varpi(t)\|^{2},  \\
	\partial_{t} \hat{v}(x, t) =&\mu \partial_{x}\hat{v}(x, t) +\hat{c}(x,t)u(x,t)  \nonumber \\
	&+ \rho \epsilon(x, t)\|\varpi(t)\|^{2}, \\
	\hat{u}(0, t)  =& r\hat{v}(0, t), \\
	\hat{v}(1, t)  =& U(t), \label{eq:uv_identifier_U} 
\end{align}
for some design gain $\rho>0$, and where 
\begin{equation}
	e (x,t)= u(x,t) - \hat{u}(x,t) , \epsilon (x,t)= v(x,t) - \hat{v}(x,t),
	\label{eq:e_identifier_def}
\end{equation}
are errors between $u$ and $v$ and their estimates $\hat{u}$ and $\hat{v}$
, $\hat{c}$ is estimate $c$. We define $\varpi(x, t)=[u(x, t),v(x, t)]^{T},$
for some initial conditions
$\hat{u}_{0}, \hat{v}_{0} \in L^2([0,1]).$
The error signals \eqref{eq:e_identifier_def} can straightforwardly be shown to have dynamics
\begin{align}
	\partial_{t} \hat{e}(x, t) =&-\lambda \partial_{x}\hat{e}(x, t)-\rho e(x, t)\|\varpi(t)\|^{2},  \label{eq:e_1}\\
	\partial_{t} \hat{\epsilon}(x, t) =&\mu \partial_{x}\hat{\epsilon}(x, t) +\tilde{c}(x,t)u(x,t)  \nonumber \\
	&- \rho \epsilon(x, t)\|\varpi(t)\|^{2}, \\
	e(0, t)  =&r\epsilon(0,t), \label{eq:e_3}\\
	\epsilon(1, t) =&0, \label{eq:e_4}
\end{align}
where $\tilde{c}(x,t)=c(x,t)-\hat{c}(x,t)$. We choose the following update laws
\begin{equation}\label{eq:lawc12}
	\hat{c}_{t}(x,t)  =\operatorname{Proj}_{\bar{c}}\left\{\gamma_1  e^{\gamma x}\varepsilon(x, t)u(x, t) , \hat{c}(x,t)\right\}, \\
\end{equation}
where $\gamma,\gamma_1 >0$ are scalar design gains. We consider the following adaptive backstepping transformation
\begin{align}\label{eq:backstepping_1}
	w(x,t)=&\hat{u}(x,t),\\
	z(x,t)=&\hat{v}(x,t)-\int_{0}^{x}\breve{K}^{u}(x,\xi,t)\hat{u}(\xi,t)d\xi \nonumber\\
	&- \int_{0}^{x}\breve{K}^{v}(x,\xi,t)\hat{v}(\xi,t)d\xi =: T[\hat{u},\hat{v}](x,t), \label{eq:backstepping_2}
\end{align}
where the kernels $\breve{K}^{u}$ and $\breve{K}^{v}$ satisfy the following kernel functions
\begin{align}\label{eq:kernel_1}
	\mu \breve{K}_{x}^{u}(x, \xi, t)= &\lambda \breve{K}_{\xi}^{u}(x, \xi, t)+\hat{c}(\xi,t)\breve{K}^{v}(x, \xi, t),\\
	\mu\breve{K}_{x}^{v}(x, \xi, t)= &- \mu \breve{K}_{\xi}^{v}(x, \xi, t), \\
	\breve{K}^{u}(x, x, t)= & -\frac{\hat{c}(x,t)}{\lambda+\mu}, \\
	\breve{K}^{v}(x, 0, t)= & \dfrac{\lambda {r}}{\mu}\breve{K}^{u}(x, 0, t). \label{eq:kernel_4}
\end{align}
The kernels are defined over the triangular domain $\mathcal{T} = \{ (x, \xi) \; | \; 0 \leq \xi \leq x \leq 1 \}.$ The inverse transformation of \eqref{eq:backstepping_2} is as follows 
\begin{align}
	\hat{v}(x,t)=&z(x,t)+\int_{0}^{x}\breve{L}^u(x,\xi,t)w(\xi,t)d \xi\nonumber\\
	&+\int_{0}^{x}\breve{L}^v(x,\xi,t)z(\xi,t)d\xi=:T^{-1}[w,z],\label{eq:inverseTrans}
\end{align}
where $T^{-1}$ is an operator similar to $T$. According to the paper \citep{valenzuela2006boundary}, there exits a constant $\bar{L}$ so that $\|\breve{L}^u\|_{\infty}\leq \bar{L}$, $\|\breve{L}^v\|_{\infty}\leq \bar{L}$.
The nominal adaptive controller is designed as 
\begin{equation}\label{U_exact}
	U(t)=\int_{0}^{1}(\breve{K}^u(1,\xi,t)\hat{u}(\xi,t) +\breve{K}^v(1,\xi,t)\hat{v}(\xi,t)d\xi .
\end{equation} Using the transformation \eqref{eq:backstepping_1} and \eqref{eq:backstepping_2}, we get the following target system
\begin{align}
	w_t(x,t) =&-\lambda w_x(x,t)+ \rho e(x,t)\|\varpi(t)\|^2, \label{eq:target_1}\\
	z_{t}(x, t)= & \mu z_{x}(x, t)-\int_{0}^{x} \breve{K}_{t}^{u}(x, \xi, t) w(\xi, t) d \xi \nonumber\\
	&-\int_{0}^{x} \breve{K}_{t}^{v}(x, \xi, t) T^{-1}[w, z](\xi, t) d \xi \nonumber\\
	& +T\left[0, \hat{c}e \right](x, t)+\rho T[e, \epsilon](x, t)\|\varpi(t)\|^{2}, \label{eq:target_2}\\
	w(0, t)= & r z(0, t), \label{eq:target_3}\\
	z(1, t)= & 0.  \label{eq:target_4}
\end{align}

\section{Neural Operator for Backstepping Kernels}
The neural operator is employed for approximating the function mapping. In the section, we introduce DeepONet for the mapping from parameter $\hat{c}$ to kernels $\breve{K}^u,\breve{K}^v$. 
\begin{lem}\label{lem1}\textit{(Existence and bound for $\breve{K}^u(x, y, t)$, $\breve{K}^v(x, y, t)$, $\breve{K}^w_t(x, y, t)$, $\breve{K}^z_t(x, y, t)$)}  
	Let $\|\hat{c}(\cdot, t)\|_{\infty} \leq \bar{c}$, $ \forall (x, t) \in [0, 1] \times \mathbb{R}_+$. Then, for any fixed $t\in\mathbb{R}_+$ and for any $\lambda,\mu \in \mathbb{R}_+$, $\hat{c} \in C^1([0,1])$,
	the equations \eqref{eq:kernel_1}--\eqref{eq:kernel_4} have unique $(C^1(\mathcal T))^2$ solution with the following  property
	\begin{align}
		|\breve{K}^{u}(x,\xi,t)|\leq& \bar{K}, \quad t \geq 0,\\
		|\breve{K}^{v}(x,\xi,t)|\leq& \bar{K}, \quad t \geq 0,
		\label{equ-k-bouded} \\
		\|\breve{K}^u_{t}(x,\xi,t)\| \in& \mathcal{L}_2,\quad  t \geq 0, \\
		\|\breve{K}^v_{t}(x,\xi,t)\| \in& \mathcal{L}_2, \quad  t \geq 0,    
	\end{align}
	where $\bar{K}>0$ is positive constant.
\end{lem}

\begin{pf}\rm
	The proof can be found in  \citep{di2013stabilization}.
\end{pf}

\begin{thm}\label{TheoDeepONet0}
	[Universal approximation theorem of DeepONet \citep{Approximation}].For sets \( X \subset \mathbb{R}^{d_x} \) and \( Y \subset \mathbb{R}^{d_y} \), they are compact sets of vectors \( x \in X \) and \( y \in Y \), respectively. Define \( \mathcal{U}: X \mapsto U \subset \mathbb{R}^{d_u} \) and \( \mathcal{V}: Y \mapsto V \subset \mathbb{R}^{d_v} \) as sets of continuous functions \( u(x) \) and \( v(y) \), respectively, and assume that \( \mathcal{U} \) is also a compact set. If the operator \( \mathcal{G}: \mathcal{U} \mapsto \mathcal{V} \) is continuous, then for any \( \mathcal{E}  > 0 \), there exist \( a^* \) and \( b^* \in \mathbb{N} \) such that when \( a \geq a^*, b \geq b^* \), there exist neural networks \( f^{\mathcal{N}}(\cdot; \theta^{(i)}), g^{\mathcal{N}}(\cdot; \theta^{(i)}), i=1, \dots, b \) with parameters \( \theta^{(i)}, \vartheta^{(i)} \), and their corresponding \( x_j \in X, j=1, \dots, a \), such that \( \mathbf{u}_a = (u(x_1), u(x_2), \dots, u(x_a))^T \) satisfies the following conditions.
	\begin{equation}
		|\mathcal{G}(u)(y)-\mathcal{G}_{\mathbb{N}}(\mathbf{u}_a)(y)|\leq \mathcal{E} ,
	\end{equation}
	where 
	\begin{equation}
		\mathcal{G}_{\mathbb{N}}(\mathbf{u}_a)(y)=\sum\limits_{i=1}^{b} g^{\mathcal{N}}(\mathbf{u};\vartheta^{(i)})f^{\mathcal{N}}(y;\theta^{(i)}),
	\end{equation}
	for all $u\in \mathcal{U}$ and $y\in Y$ of $\mathcal{G}(u)\in \mathcal{V}$.
\end{thm}
\par Next, we define the operator $\mathcal{K}: C^1([0, 1]) \mapsto (C^1(\mathcal{T}))^2$
\begin{align}\label{eq:operatorK}
	\mathcal{K}(\hat{c}) (\cdot,t):=(\breve{K}^u(x,\xi),\breve{K}^v(x,\xi)).
\end{align}

\begin{thm}\label{th:NO}
	{\em [Existence of a NO to approximating the kernel]} Consider the NO defined in \eqref{eq:operatorK}, along with \eqref{eq:kernel_1}-\eqref{eq:kernel_4}.
	Fix $t>0$. Then for all $\mathcal{E}  > 0$ there exists a NO $\mathcal{\hat{K}}:C^1([0, 1]) \mapsto (C^1(\mathcal{T}))^2$ such that for all $(x,\xi) \in \mathcal{T}$,
	\begin{align*}
		|\mathcal{K}(\hat{c})(\cdot,t)-\hat{\mathcal{K}}(\hat{c})(\cdot,t)|\leq\mathcal{E}.
	\end{align*}
	\end{thm}
	\begin{pf}
		The continuity of the operator \(\mathcal{K}\) is derived from Lemma \ref{lem1}. And this result is based on Theorem 2.1 proposed by B. Deng et al. in their study \citep{Approximation}.
	\end{pf}

	\section{ Stabilization under DeepONet-Approximated Gain Feedback}\label{section_DeepONet}\label{sec_DeepONet}
	We will establish the system stability property of the adaptive backstepping controller using the approximate estimated kernels $\hat{K}^u$ and $\hat{K}^v$ in the following theorem.
	
	\begin{thm}
		{\em [Stabilization under approximate adaptive backstepping control.]} Consider the system \eqref{eq:uv_def}-\eqref{eq:uv_def_v0}, there exists a $\mathcal{E} _0$ such that for all NO approximations $\hat{K}^u,\hat{K}^v$ of accuracy $\mathcal{E}  \in (0,\mathcal{E}_0)$ provided by Theorem \ref{th:NO}, the plant \eqref{eq:uv_def}-\eqref{eq:uv_def_v0} in feedback with the adaptive control law 
		\begin{equation}\label{eq:U_NO}
			U(t)=\int_{0}^{1}\hat{K}^u(1,\xi,t)\hat{u}(\xi,t) +\hat{K}^v(1,\xi,t)\hat{v}(\xi,t)d\xi ,
		\end{equation} 
		along with the update law for $\hat{c}$ given by \eqref{eq:lawc12} with any Lipschitz initial condition $\hat{c}_{0}=\hat{c}(\cdot,0)$ such that $\hat{c}_{0}\leq \bar{c}$ and the passive identifier $\hat{u}$, $\hat{v}$ given by \eqref{eq:uv_identifier_u}-\eqref{eq:uv_identifier_U} with any initial condition $\hat{u}_0=\hat{u}(\cdot,0),\hat{v}_0=\hat{v}(\cdot,0)$ such that $\|\hat{u}_0\|<\infty,\|\hat{v}_0\|<\infty$, the following properties hold:
		\begin{align}
			&\|u\|,\|v\|,\|\hat{u}\|,\|\hat{v}\|, \|{u}\|_{\infty},\|{v}\|_{\infty}, \|\hat{u}\|_{\infty},\|\hat{v}\|_{\infty}, \in L_2\cap L_{\infty},\\
			& \|{u}\|_{\infty},\|{v}\|_{\infty},\|\hat{u}\|_{\infty},\|\hat{v}\|_{\infty} \mapsto 0.
		\end{align}
		Moreover, for the equilibrium $(u,v,\hat{u},\hat{v},\hat{c})=(0,0,0,0,{c})$ the following global stability estimate holds 
		\begin{equation}\label{eq:S(t)}
			S(t)\leq 2\dfrac{k_2}{k_1}\theta_2 S(0) e^{\theta_1 k_2S(0)}, t>0,
		\end{equation}
		where
		\begin{equation}
			S(t) := \|u\|^2+\|v\|^2+\|\hat{u}\|^2+\|\hat{v}\|^2+\|\tilde{c}\|^2,
		\end{equation}
		and $k_1$, $k_2$, $\theta_1$ and $\theta_2$ are strictly positive constants. 
	\end{thm}
	\begin{pf}
		\par We consider the adaptive backstepping transformations  \eqref{eq:backstepping_1} and \eqref{eq:backstepping_2}, which lead to the following target system 
		\begin{align}
			w_t(x,t) =&-\lambda w_x(x,t)+ \rho e(x,t)\|\varpi(t)\|^2, \label{eq:target_NO_1}\\
			z_{t}(x, t)= & \mu z_{x}(x, t)-\int_{0}^{x} \breve{K}_{t}^{u}(x, \xi, t) w(\xi, t) d \xi \nonumber\\
			&-\int_{0}^{x} \breve{K}_{t}^{v}(x, \xi, t) T^{-1}[w, z](\xi, t) d \xi \nonumber\\
			& +T\left[0, \hat{c}e \right](x, t)+\rho T[e, \epsilon](x, t)\|\varpi(t)\|^{2}, \label{eq:target_NO_2}\\
			w(0, t)= & r z(0, t),\label{eq:target_NO_3} \\
			z(1, t)= & \int_{0}^{1}\tilde{K}^u(1,\xi,t)\hat{u}(\xi,t)d\xi \nonumber\\
			&+ \int_{0}^{1}\tilde{K}^v(1,\xi,t)\hat{v}(\xi,t)d\xi:=\Gamma(t), \label{eq:target_NO_4}
		\end{align}
		We use the following Lyapunov function candidate
		\begin{align}
			V(t):=V_1 + a V_2
		\end{align}
		where $a>0$ and 
		\begin{align}
			V_1(t)&=\int_{0}^{1}e^{-\delta x}w^2(x,t)dx, \label{eq:V_1}\\
			V_2(t)&=\int_{0}^{1}e^{kx}z^2(x,t)dx. \label{eq:V_2}
		\end{align}
		$\delta,k$ are constants to be decided, assuming $\delta \geq 1$. 
		\par (a) Bound on $\dot{V}_1$
		\par From differentiating \eqref{eq:V_1} with respect to time, inserting the dynamics \eqref{eq:target_1}, and integrating by parts, we find
		\begin{align}
			\dot{V}_{1}(t)= & -\lambda e^{-\delta} w^{2}(1, t)+\lambda w^{2}(0, t)-\lambda \delta \int_{0}^{1} e^{-\delta x} w^{2}d x \nonumber \\
			& +2 \int_{0}^{1} e^{-\delta x} w(x, t) \rho e(x, t)\|\varpi(t)\|^{2} d x.
		\end{align}
		Consider the last term, we use Young’s and  Minkowski’s inequalities, then get
		\begin{align}
			\dot{V}_{1}(t)\leq & -\lambda e^{-\delta} w^{2}(1, t)+\lambda w^{2}(0, t)-\lambda \delta \int_{0}^{1} e^{-\delta x} w^{2}d x \nonumber \\
			& +\rho_{1} \rho^{2}\|w(t)\|^{2}\|e(t)\|^{2}\|\varpi(t)\|^{2} \nonumber\\
			& +\frac{4}{\rho_{1}}\left(\left(1+2 A_{3}^{2}\right)\|w(t)\|^{2}\right)\nonumber\\
			&+\frac{4}{\rho_{1}}\left(2 A_{4}^{2}\|z(t)\|^{2}+\|e(t)\|^{2}+\|\epsilon(t)\|^{2}\right)
		\end{align}
		for some arbitrary $\rho_1,A_3,A_4 > 0$. Choosing $\rho_1 = e^\delta$, we find
		\begin{align}
			\dot{V}_1(t) \leq & \lambda r^2 z^{2}(0, t)-\left[\lambda \delta-4-8A_3^2\right] V_1(t)\nonumber\\
			&+h_{3} V_{5}(t)+l_{1}(t) V_1(t)+l_{2}(t).
		\end{align}
		where
		\begin{align}
			l_1(t)=&e^{2\delta}\rho^2\|e(t)\|^{2}\|\varpi(t)\|^{2},\\
			l_2(t)=&4e^{-\delta}\rho^2(\|e(t)\|^{2}+\|\epsilon(t)\|^{2})
		\end{align}
		\par (b) Bound on $\dot{V}_2$
		\par From differentiating \eqref{eq:V_2} with respect to time, inserting the dynamics \eqref{eq:target_2}, and integrating by parts, we find
		\begin{align}
			\dot{V}_2(t)= & \mu e^{k} z^{2}(1, t)-\mu z^{2}(0, t)-\mu k \int_{0}^{1} e^{k x} z^{2}(x, t) d x \nonumber\\
			& -2 \int_{0}^{1} e^{k x} z(x, t) \int_{0}^{x} \hat{K}_{t}^{u}(x, \xi, t) w(\xi, t) d \xi d x \nonumber\\
			& -2 \int_{0}^{1} e^{k x} z(x, t) \int_{0}^{x} \hat{K}_{t}^{v}(x, \xi, t) T^{-1}[w, z]d \xi d x \nonumber\\
			& +2 \int_{0}^{1} e^{k x} z(x, t) T\left[0, \hat{c}e\right](x, t) d x \nonumber\\
			& +2 \rho \int_{0}^{1} e^{k x} z(x, t) T[e, \epsilon](x, t)\|\varpi(t)\|^{2} d x .
		\end{align}
		Using Young’s inequality
		\begin{align}
			\dot{V}_{5}(t) \leq & \mu e^{k} z^{2}(1, t)-\mu z^{2}(0, t)-\left[k \mu-3\right] \int_{0}^{1} e^{k x} z^{2}d x \nonumber\\
			& +\int_{0}^{1} e^{k x}\left(\int_{0}^{x} \hat{K}_{t}^{u}(x, \xi, t) w(\xi, t) d \xi\right)^{2} d x \nonumber\\
			& +2 \int_{0}^{1} e^{k x}\left(\int_{0}^{x} \hat{K}_{t}^{v}(x, \xi, t) T^{-1}[w, z](\xi, t) d \xi\right)^{2} d x \nonumber\\
			& +\int_{0}^{1} e^{k x} T^{2}\left[0, \hat{c} e\right](x, t) d x \nonumber\\
			& +\rho^{2} e^{\delta+2 k} \int_{0}^{1}(z(x, t) T[e, \epsilon](x, t))^{2}\|\varpi(t)\|^{2} d x \nonumber\\
			& +e^{-\delta}(\|u(t)\|+\|v(t)\|)^{2}
		\end{align}
		and Cauchy–Schwarz’ and Minkowski’s inequalities
		\begin{align}
			\dot{V}_2(t) \leq & -\mu z^{2}(0, t)+4(1+A_3)^{2} V_1-\left[k \mu-3-4A_4^2\right] V_2 \nonumber\\
			& +l_{3}(t) V_1+l_{4}(t) V_2+l_{5}(t)+\mu e^k\Gamma(t)^2
		\end{align}
		where $A_1,A_2>0$,
		\begin{align}
			l_3(t)=&e^{k+\delta}\|\hat{K}_t^u\|^2+2e^{k+\delta}\|\hat{K}_t^v\|^2,\\
			l_4(t)=&2\rho^2e^{\delta+2k}(A_1^2\|e(t)\|^2+A_2^2\|\epsilon(t)\|^2)\|\varpi(t)\|^{2}\nonumber\\
			&+2e^{k}\|\hat{K}_t^v\|^2A_4^2,\\
			l_5(t)=&2e^kA_2^2\bar{c}\|e(t)\|^{2}+4e^{-\delta}(\|e(t)\|^{2}+\|\epsilon(t)\|^{2}).
		\end{align}
		According to \eqref{eq:inverseTrans}, we have
		$
		\Gamma(t)^2 \leq 2\mathcal{E}^2 L_1V(t),$
		where $L_1=2\bar{L}^2+3\bar{L}+1$. Thus, we find
		\begin{align}\label{eq:V6_dot1_1}
			\dot{V}(t) \leq &-[a\mu-\lambda r^2]z^2(0,t)+(l_1(t)+al_3(t))V_1(t)\nonumber\\
			&-(\lambda\delta-4-8A_3^2-4a(1+A_3)^2)V_1(t)\nonumber\\
			&-(ak\mu-a(3+4A_4^2)-8e^{-\delta}A_4^2)V_2(t)+l_2(t)\nonumber\\
			&+al_5(t)+2\mu e^k \mathcal{E}^2L_1V(t),
		\end{align}
		Let 
		$$a=\dfrac{\lambda r^2+1}{\mu}, \delta>\max \left\{1,\dfrac{4+8A_3^2+4a(1+A_3)^2}{\lambda}\right\},$$ $$k>\dfrac{8e^{-\delta}A_4^2+a(3+4A_4^2)}{a\mu}.$$ Thus, we obtain the following upper bound calculation 
		\begin{equation}\label{eq:V6_dot1_1}
			\dot{V}(t) \leq-\left[d-2\mu e^k\mathcal{E}^2 L_1 \right] V(t)+l_6(t) V(t)+l_{7}(t),
		\end{equation}
		for positive constant $d$ and the nonnegative, integrable functions $l_{6}(t)$ and $l_7(t)$
		\begin{align}
			l_{6}(t) = &  \mathop{\max}
			\left\{l_1(t)+al_3(t),l_4(t)\right\}, \label{eq:l6}\\
			l_7(t)=&-z^2(0,t)+l_2(t)+al_5(t),\label{eq:l7}
		\end{align}
		We introduce 
		\begin{equation}
			\mathcal{E} _0:=\dfrac{\sqrt{2d-1}}{2\sqrt{\mu e^kL_1}}.
		\end{equation} 
		
		Thus, if we choose $\mathcal{E} \in (0,\mathcal{E}_0)$ we have $d-2\mu e^k\mathcal{E}^2 L_1 >{1}/{2}>0$.
		It then follows from Lemma B.6 \citep{krstic1995nonlinear} that
		$$
		V \in \mathcal{L}_1 \cap \mathcal{L}_{\infty},
		$$
		and hence
		$
		\|w\|,\|z\| \in \mathcal{L}_2 \cap \mathcal{L}_{\infty} .
		$
		Due to the invertibility of the backstepping transformation                                            
		$$
		\|\hat{u}\|,\|\hat{v}\| \in \mathcal{L}_2 \cap \mathcal{L}_{\infty} .
		$$
		From $\|{e}\|,\|{\epsilon}\| \in \mathcal{L}_2 \cap \mathcal{L}_{\infty}$\citep{anfinsen2019adaptive} , it follows that 
		$$
		\|u\|,\|v\| \in \mathcal{L}_2 \cap \mathcal{L}_{\infty} .
		$$
		\par Then, we proceed by showing pointwise boundedness, square integrability and convergence to zeros. The paper \citep{vazquez2011backstepping} proved that the system \eqref{eq:uv_def}-\eqref{eq:uv_def_v0} is equivalent to the following system through an invertible backstepping transformation.
		\begin{align}\label{alpha_1}
			\alpha_t(x, t)&= -\lambda \alpha_x(x, t)+h(x) \beta(0, t), \\
			\beta_t(x, t)&=\mu \beta_x(x, t), \\
			\alpha(0, t)&= r \beta(0, t), \\
			\beta(1, t)
			&=U(t) -\int_{0}^{1}G_1(\xi)u(\xi)d\xi -\int_{0}^{1}G_2(\xi)v(\xi)d\xi, \label{alpha_4}
		\end{align}
		for some bounded functions $h,G_1,G_2\in \mathcal{C}$ fo the unknown parameters. Equations \eqref{alpha_1}-\eqref{alpha_4} can be explicitly be solved for $t>\lambda^{-1}+\mu^{-1}$ to yield 
		\begin{align}\label{eq:alpha_solution}
			\alpha(x, t)=&r \beta\left(1, t-\mu^{-1}-\lambda^{-1} x\right) \nonumber \\
			&+\lambda^{-1} \int_{0}^{x} h(\tau) \beta\left(1, t-\mu^{-1}-\lambda^{-1}(x-\tau)\right) d \tau, \\
			\beta(x, t)=&\beta\left(1, t-\mu^{-1}(1-x)\right) . \label{eq:alpha_solution_2}
		\end{align}
		From \eqref{alpha_4}, the control law $U(t)$ and $\|u\|,\|v\|,\|\hat{u}\|,\|\hat{v}\| \in$ $\mathcal{L}_2 \cap \mathcal{L}_{\infty}$, it follows that $\beta(1, \cdot) \in$ $\mathcal{L}_2 \cap \mathcal{L}_{\infty}$. Since $\beta$ and $\alpha$ are simple, cascaded transport equations, this implies
		\begin{align}
			\|\alpha\|_{\infty},\|\beta\|_{\infty} \in \mathcal{L}_2 \cap \mathcal{L}_{\infty}, \quad\|\alpha\|_{\infty},\|\beta\|_{\infty} \rightarrow 0.
		\end{align}
		With the invertibility of the transformation, then yields
		\begin{equation}
			\|u\|_{\infty},\|v\|_{\infty} \in \mathcal{L}_2 \cap \mathcal{L}_{\infty}, \quad\|u\|_{\infty},\|v\|_{\infty} \rightarrow 0.
		\end{equation}

		From the structure of the identifier \eqref{eq:uv_def}, we will also have  $\hat{u}(x, \cdot), \hat{v}(x, \cdot) \in \mathcal{L}_{\infty} \cap \mathcal{L}_{2}$ , and hence
		\begin{align}
			\|\hat{u}\|_{\infty}, \|\hat{v}\|_{\infty} \in \mathcal{L}_{\infty} \cap \mathcal{L}_{2}, 
			\quad \|\hat{u}\|_{\infty},\|\hat{v}\|_{\infty} \mapsto 0.
		\end{align}

		\par Next, we consider the global stability. We will use the Lyapunov function to show that the system's state remains stable over time,
		\begin{align}\label{eq:V_1_1}
			V_3(t) =\int_{0}^{1}e^{-\gamma x}e^2(x,t)dx + \int_{0}^{1}e^{\gamma x}\varepsilon^2(x,t)dx + \dfrac{\|\tilde{c}(t)\|^2}{\gamma_1},
		\end{align}
		Inserting the dynamics \eqref{eq:e_1}-\eqref{eq:e_4}, integrating by parts, and using the property $-\tilde{c} \operatorname{Proj}_{\bar{c}}\left\{\gamma_1  e^{\gamma x}\varepsilon(x, t)u(x, t) , \hat{c}(x,t)\right\} \leq -\tilde{c}\gamma_1e^{\gamma x}\epsilon u$(\citep{anfinsen2019adaptive}), we find it leads to the following upper bound:
		\begin{align}
			\dot{V}_3(t) \leq & -\lambda e^{-\gamma} e^{2}(1, t)-(\mu-\lambda r^2 )e^{2}(0, t)\nonumber\\
			&-\gamma \lambda e^{-\gamma}\|e(t)\|^{2}-2\rho e^{-\gamma}\|e(t)\|^{2}\|\varpi(t)\|^{2}\nonumber\\
			&-\mu \gamma\|\epsilon(t)\|^{2}-2 \rho e^{\gamma}\|\epsilon(t)\|^{2}\|\varpi(t)\|^{2},\label{eq:dotV_1}
		\end{align}
		which shows that ${V}_3(t)$ is non-increasing and hence bounded. Thus implies that the ${V}_3(t)<V_3(0)$ and limit $\lim_{t\to \infty}V_3(t)=V_{3,\infty}$ exists. By integrating \eqref{eq:dotV_1} from zero to infinity, we obtain the following upper bound:
		\begin{align}
			&\lambda e^{-\gamma} \int_{0}^{\infty}e^{2}(1, \tau)d\tau+(\mu-\lambda r^2) \int_{0}^{\infty}e^{2}(0, \tau) d\tau \nonumber\\
			&+\gamma \lambda  e^{-\gamma}\int_{0}^{\infty}\|e(\tau)\|^{2}d\tau +2 \rho e^{-\gamma}\int_{0}^{\infty}\|e(\tau)\|^{2}\|\varpi(\tau)\|^{2}d\tau\nonumber\\
			& +\mu \gamma \int_{0}^{\infty}\|\epsilon(\tau)\|^{2}d\tau+2 \rho e^{\gamma}\int_{0}^{\infty}\|\epsilon(\tau)\|^{2}\|\varpi(\tau)\|^{2}d\tau \nonumber\\
			&\leq V_3(0).
			\label{eq:dotV_1(0)}
		\end{align}
		
		From \eqref{eq:l6} and \eqref{eq:l7}, it can be concluded that there are constants $\theta_1 >0 $ and  $\theta_2>1$ such that
		\begin{align}\label{eq:l6_1}
			\|l_6\|_1 \leq& \theta_1 V_3(0), \\
			\|l_{10}\|_1 \leq& \theta_2 V_3(0). \label{eq:l10_1}
		\end{align}
		Recalling \eqref{eq:V6_dot1_1}, we have that 
		\begin{align}
			\dot{V}(t) \leq -\frac{1}{2}V(t) +l_6(t)V(t) +l_{7}(t).
		\end{align}
		We also have that from Lemma B.6 in the paper\citep{krstic1995nonlinear} 
		\begin{equation}\label{eq:V leq}
			V(t)\leq (e^{-\frac{1}{2}t}V(0)+\|l_{7}\|_1)e^{\|l_6\|_1}.
		\end{equation}
		\par We then introduce the function
		\begin{equation}\label{eq:V3}
			V_4(t):=V_3(t)+V(t).
		\end{equation}
		Noticing that
		\begin{equation}\label{eq:V1 leq}
			V_3(t)\leq 	V_3(0)\leq \theta_2 V_3(0)e^{\theta_1 V_3(0)},		
		\end{equation}
		we achieve from \eqref{eq:V leq}, \eqref{eq:V1 leq}, \eqref{eq:l6_1} and \eqref{eq:l10_1} the following
		\begin{align}\label{eq:V3 leq}
			V_4(t)= & V(t)+V_3(t)\nonumber \\
			\leq &  (e^{-\frac{1}{2}t}V(0)+\|l_{7}\|_1)e^{\|l_6\|_1} +\theta_2 V_3(0)e^{\theta_1 V_3(0)}	\nonumber	\\
			\leq & (\theta_2 V(0)+\theta_2V_3(0))e^{\|l_6\|_1} +\theta_2 V_3(0)e^{\theta_1 V_3(0)} \nonumber\\
			\leq &2\theta_2 V_4(0) e^{\theta_1 V_4(0)}.
		\end{align}
		This Lyapunov functional can be represented by an equivalent norm, and the bounds of this equivalent norm are determined by two positive constants $k_1 > 0$ and $k_2 > 0$,
		\begin{equation}
			k_1 S(t)\leq V_4(t) \leq k_2 S(t).
		\end{equation}
		So we have
		\begin{equation}
			S(t)\leq 2\dfrac{k_2}{k_1}\theta_2 S(0) e^{\theta_1 k_2S(0)}.
		\end{equation}
		
	\end{pf}
	\section{Simulations}\label{sec_Simulation}
	This section we analyze the performance of the proposed NO-based adaptive control law for the ARZ traffic PDE system through simulations on a L=600m road over T=300 seconds. The parameters are set as follows: free-flow velocity \(v_m = 40\) m/s, maximum density \(\rho_m = 160\) veh/km, equilibrium density \(\rho_1^* = 120\) veh/km, driver reaction time \(\tau = 60\) s. Initial conditions are sinusoidal inputs \(\rho_1(x,0) = \rho_1^* + 0.1\sin(\frac{3\pi x}{L})\rho_1^*\) and \(v_1(x,0) = v_1^* - 0.01\sin(\frac{3\pi x}{L})v^*\) to mimic stop-and-go traffic. Let $\gamma_0=1$. 
	\par
\begin{figure}
		\begin{center}
			\includegraphics[width=7cm]{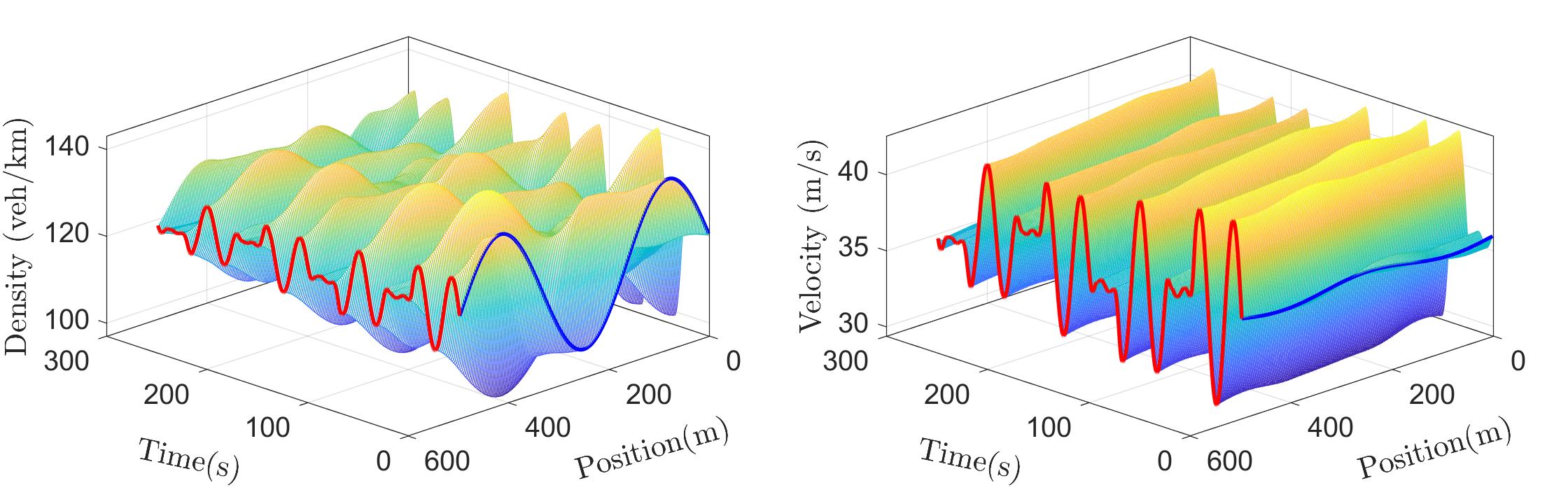}    
			\caption{Density and velocity evolution of open-loop ARZ traffic system }
			\label{ARZ_openloop}
		\end{center}
	\end{figure}

\begin{figure}[ht] 
		\centering
		\begin{subfigure}{0.15\textwidth}
			\includegraphics[width=\textwidth]{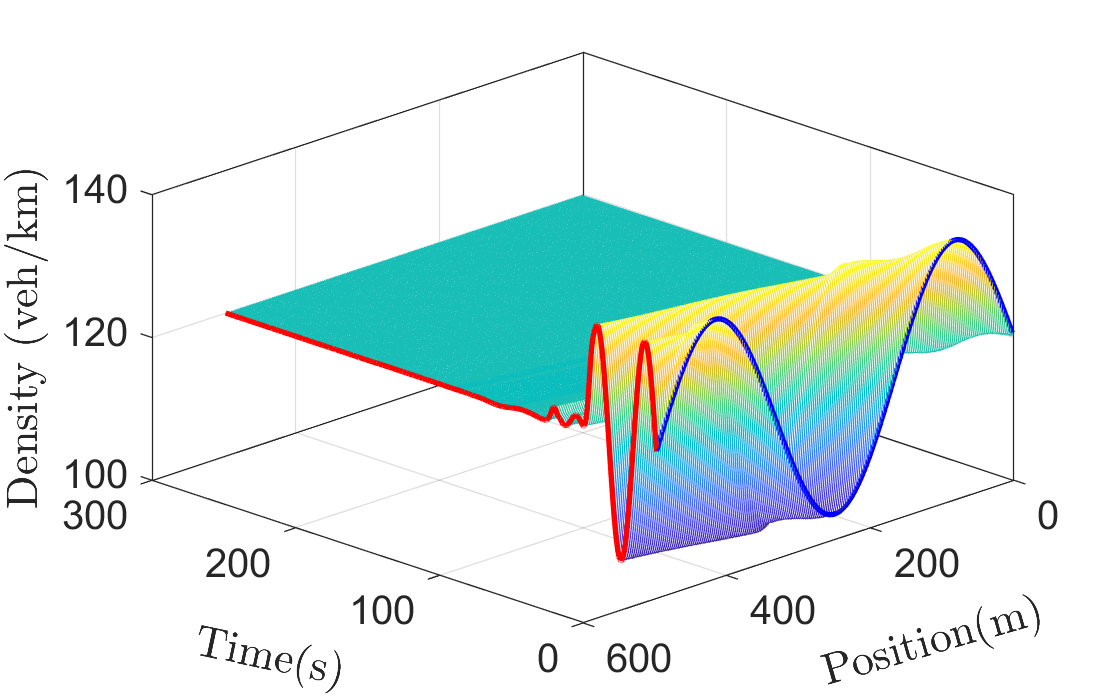} 
			\caption{$\rho_1(x,t)$}
			\label{fig:rho1}
		\end{subfigure}
		\hfill
		\begin{subfigure}{0.15\textwidth}
			\includegraphics[width=\textwidth]{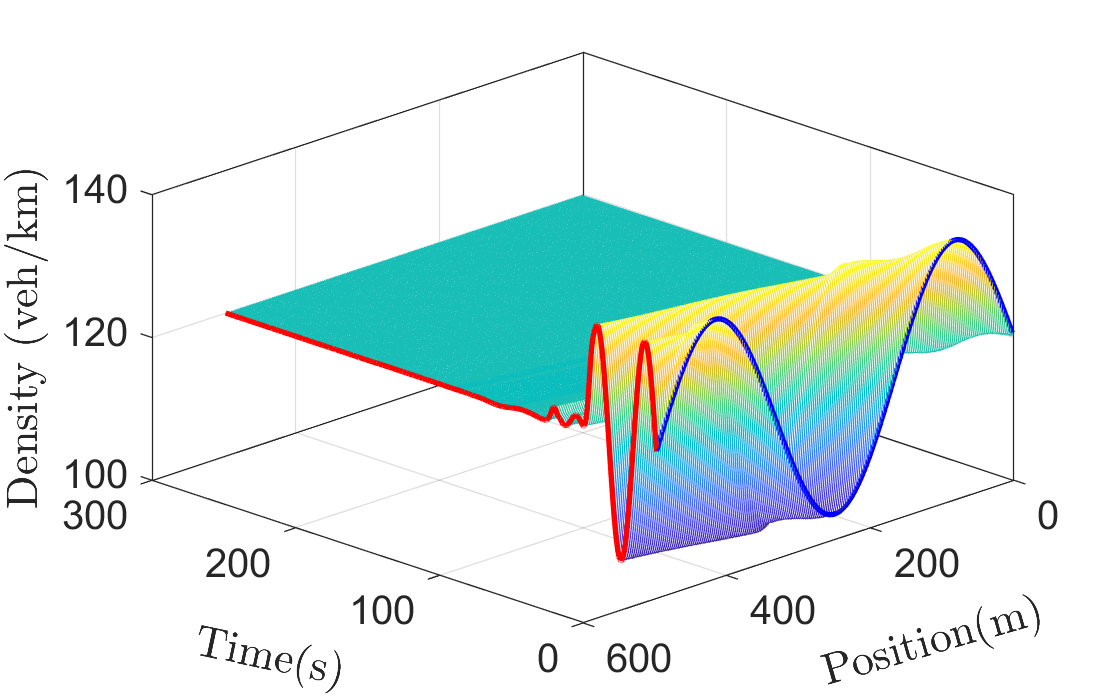}
			\caption{$\rho_{1NO}(x,t)$}
			\label{fig:rho1NO}
		\end{subfigure}
		\hfill
		\begin{subfigure}{0.15\textwidth}
			\includegraphics[width=\textwidth]{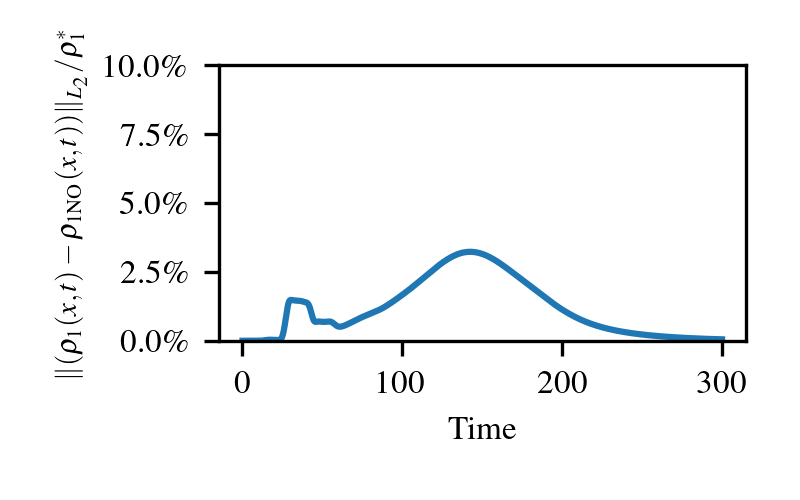}
			\caption{Relative $L_2$ error}
			\label{fig:rho_error}
		\end{subfigure}

		\vspace{1em} 
		
		\begin{subfigure}{0.15\textwidth}
			\includegraphics[width=\textwidth]{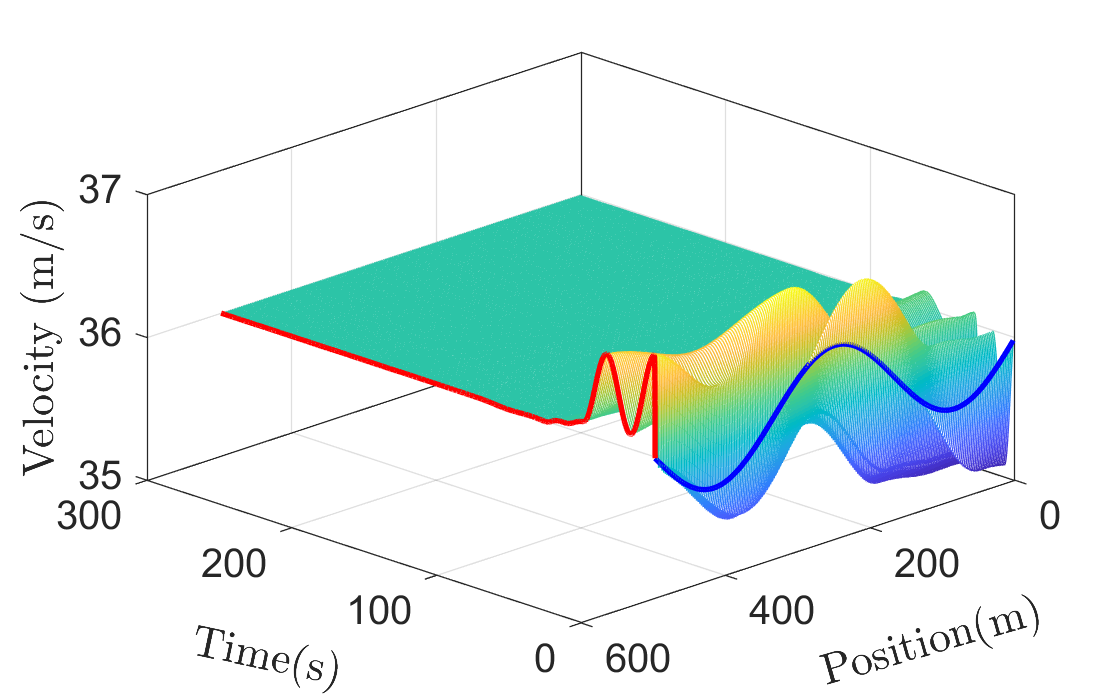}
			\caption{$v_1(x,t)$}
			\label{fig:v1}
		\end{subfigure}
		\hfill
		\begin{subfigure}{0.15\textwidth}
			\includegraphics[width=\textwidth]{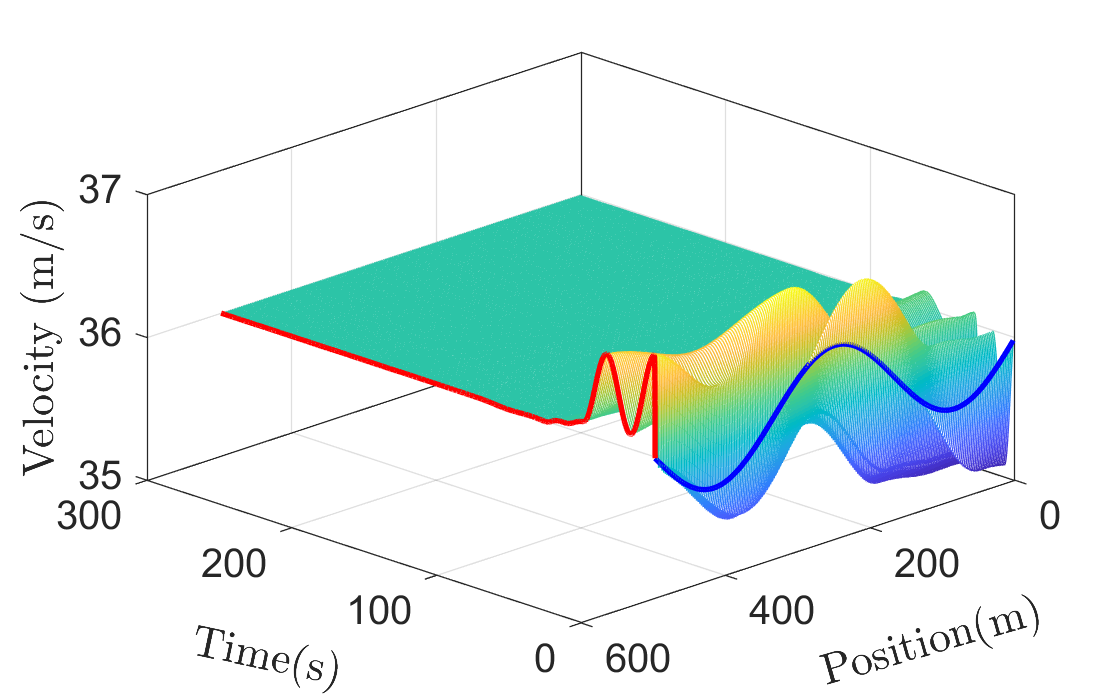}
			\caption{$v_{1NO}(x,t)$}
			\label{fig:v1NO}
		\end{subfigure}
		\hfill
		\hfill
		\begin{subfigure}{0.15\textwidth}
			\includegraphics[width=\textwidth]{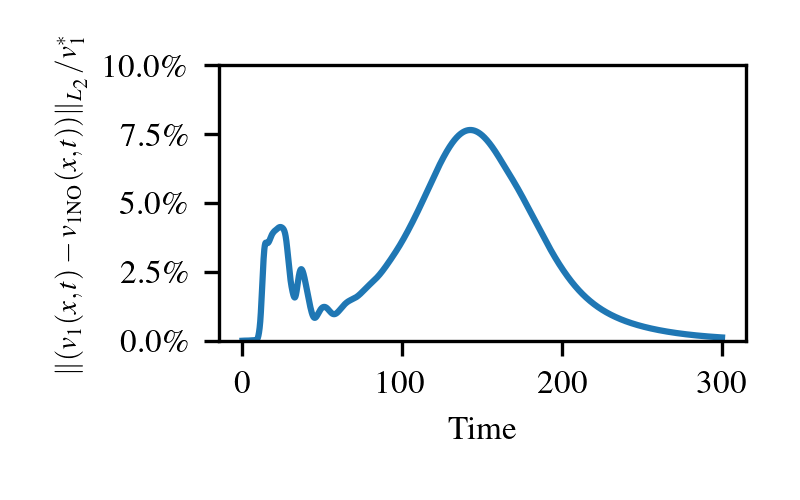}
			\caption{Relative $L_2$ error}
			\label{fig:v_error}
		\end{subfigure}
		
		\caption{Simulation results for density and velocity with feedback controller \eqref{U_exact} and \eqref{eq:U_NO}. The results show $\rho_1(x,t)$, $\rho_{1NO}(x,t)$, and their relative $L_2$ error, as well as $v_1(x,t)$, $v_{1NO}(x,t)$, and their relative $L_2$ error.}
		\label{fig:rho_v}
	\end{figure}

\begin{figure}
		\begin{center}
			\includegraphics[width=7cm]{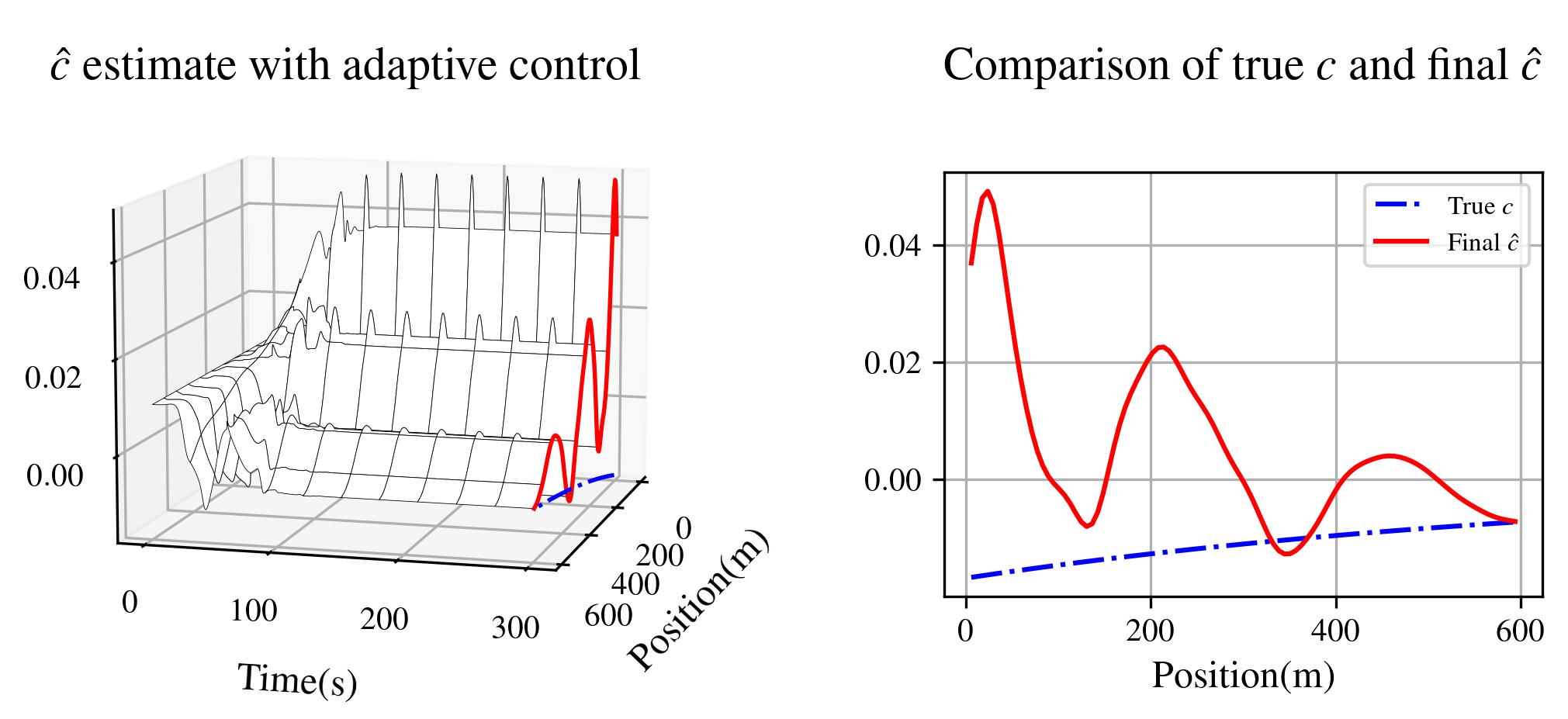}    
			\caption{Parameter estimation of $\hat{c}$ in ARZ traffic system and comparision of true $c$ and $\hat{c}$}
			\label{c_Comparison}
		\end{center}
	\end{figure}

 To generate a sufficient dataset for training, we use 10 different \(c(x)\) functions with \(\tau \in U[50,70]\) and simulate the resulting PDEs under the adaptive controller for \(T=300\) seconds. We sub-sample each \((c, \hat{K}^{u}, \hat{K}^{v})\) pair every 0.1 seconds, resulting in a total of 30,000 distinct \((c, \hat{K}^{u}, \hat{K}^{v})\) pairs for training the NO. 
	Figure \ref{ARZ_openloop} shows the ARZ system is open-loop unstable. Figure \ref{fig:rho_v} show the density and velocity of ARZ traffic system with the adaptive controller and NO-based adaptive controller. 
	The results indicate that both the NO-based adaptive method and the adaptive backstepping control method effectively stabilize the transportation system. The maximum error does not exceed 10$\%$. The traffic density and velocity converge to the equilibrium values of \(\rho_1^* = 120\) veh/km and \(v_1^* = 36\) m/s, respectively.   The estimated parameter \(\hat{c}\) is shown in Figure \ref{c_Comparison}. Table \ref{tab2} summarizes the maximum and mean absolute errors for backstepping kernels and traffic states, showing high accuracy in computation of  kernels and traffic states.  Notably, the NO-based adaptive control method accelerates controller computation by up to 150 times compared to traditional adaptive control methods.

	\begin{table}[hb]
		\begin{center}
			\caption{The errors of kernels and traffic states}\label{tab2}
			\begin{tabular}{ccc}
				& Max absolute error & Mean absolute error \\ \hline
				Kernel $K^u$ & $1.24 \times 10^{-3}$ & $1.06 \times 10^{-3}$ \\
				Kernel $K^v$ & $2.48 \times 10^{-3}$ & $2.33 \times 10^{-3}$ \\
				Density (veh/km) & $1.32$ & $0.51$ \\
				Speed (km/h) & $2.35$ & $0.87$ \\ \hline
			\end{tabular}
		\end{center}
	\end{table}
	
	\section{Conclusion}
	We accelerate the computation of \(2 \times 2\) coupled Goursat-form PDEs in adaptive control designs for traffic ARZ system by using DeepONet to learn adaptive control gains. Experimental results show our method improves computational efficiency by two orders of magnitude compared to traditional solvers, enhancing the real-time applicability of adaptive control in traffic congestion. Future work will incorporate real traffic data into the training of the neural operator.
	

	\bibliography{references}

\end{document}